\documentclass{elsart}

\usepackage{amssymb}
\usepackage{graphicx}

\newcommand{\comment}[1]{}
\newcommand{\halmos}{\rule{1ex}{1.4ex}}
\newenvironment{proof}{\noindent {\em Proof}.\
}{\hspace*{\fill}$\halmos$\medskip}

\newtheorem{theorem}{Theorem}
\newtheorem{itlemma}{Lemma}[section] 
\newtheorem{itproposition}[itlemma]{Proposition}
\newtheorem{itcorollary}[itlemma]{Corollary}
\newtheorem{itremark}[itlemma]{Remark}
\newtheorem{itdefinition}[itlemma]{Definition}
\newtheorem{itexercise}[itlemma]{Exercise}
\newtheorem{itexample}[itlemma]{Example}

\newenvironment{lemma}{\begin{itlemma}\rm}{\end{itlemma}} 
\newenvironment{remark}{\begin{itremark}\rm}{\end{itremark}} 
\newenvironment{corollary}{\begin{itcorollary}\rm}{\end{itcorollary}}
\newenvironment{proposition}{\begin{itproposition}\rm}{\end{itproposition}}
\newenvironment{definition}{\begin{itdefinition}\rm}{\end{itdefinition}}
\newenvironment{exercise}{\begin{itexercise}\rm}{\end{itexercise}}
\newenvironment{example}{\begin{itexample}\rm}{\end{itexample}}
\def \beq {\begin{eqnarray}}
\def \eeq {\end{eqnarray}}
\def \beqn {\begin{eqnarray*}}
\def \eeqn {\end{eqnarray*}}

\newcommand{\text}[1]{\hbox{\rm \ #1\ \/}}
\newcommand{\be}[1]{\begin{equation}\label{#1}}
\newcommand{\ee}{\end{equation}}
\newcommand{\bl}[1]{\begin{lemma}\label{#1}}
\newcommand{\ble}[1]{\begin{lemmaex}\label{#1}}
\newcommand{\br}[1]{\begin{remark}\label{#1}}
\newcommand{\bt}[1]{\begin{theorem}\label{#1}}
\newcommand{\bd}[1]{\begin{definition}\label{#1}}
\newcommand{\bp}[1]{\begin{proposition}\label{#1}}
\newcommand{\bc}[1]{\begin{corollary}\label{#1}}
\newcommand{\bfact}[1]{\begin{fact}\label{#1}}
\newcommand{\ber}[1]{\begin{exercise}\label{#1}}
\newcommand{\bex}[1]{\begin{example}\label{#1}}
\newcommand{\bem}[1]{\begin{example}\label{#1}}  
\newcommand{\ec}{\mybox\end{corollary}}
\newcommand{\efact}{\mybox\end{fact}}
\newcommand{\eer}{\mybox\end{exercise}}
\newcommand{\eex}{\mybox\end{example}}
\newcommand{\eem}{\mybox\end{example}}
\newcommand{\el}{\mybox\end{lemma}}
\newcommand{\ele}{\mybox\end{lemmaex}}
\newcommand{\er}{\mybox\end{remark}}
\newcommand{\et}{\qed\end{theorem}}
\newcommand{\ed}{\mybox\end{definition}}
\newcommand{\ep}{\mybox\end{proposition}}
\newcommand{\epr}{\end{proof}}
\newcommand{\bpr}{\begin{proof}}

\newcommand{\ecs}{\end{corollary}}
\newcommand{\eers}{\end{exercise}}
\newcommand{\eexs}{\end{example}}
\newcommand{\eems}{\end{example}}
\newcommand{\els}{\end{lemma}}
\newcommand{\eles}{\end{lemmaex}}
\newcommand{\ers}{\end{remark}}
\newcommand{\ets}{\end{theorem}}
\newcommand{\eds}{\end{definition}}
\newcommand{\eps}{\end{proposition}}

\newcommand{\R}{{\mathbb R}}

\renewcommand{\qed}{\hfill \mbox{$\halmos$}} 
\newcommand{\mybox}{\hfill \mbox{$\Box$}} 

\def \min {{\rm min \,}}

\newcommand{\abs}[1]{\left\vert #1 \right\vert}
\newcommand{\rref}[1]{(\ref{#1})}

\newcommand{\dis}{\displaystyle}
\newcommand{\vf}{\varphi}

\newcommand{\mylimsup}[1]{\raisebox{-1ex}{$\stackrel{\overline{\displaystyle\lim}}{\scriptscriptstyle #1}$}}

\newcommand{\myliminf}[1]{\raisebox{-1ex}{$\stackrel{\underline{\displaystyle\lim}}{\scriptscriptstyle #1}$}}

\begin{document}

\runauthor{Eduardo D. Sontag and Yuan Wang}
\begin{frontmatter}
\title{A cooperative system which does not satisfy the limit set dichotomy}
\author[Rutgers]{Eduardo D. Sontag\thanksref{txe}}
\address[Rutgers]{Department of Mathematics\\
Rutgers University, New Brunswick, NJ 08903, USA}
\thanks[txe]{E-mail: {\tt sontag@math.rutgers.edu}}
\author[FAU]{Yuan Wang\thanksref{txw}}
\address[FAU]{Department of Mathematical Sciences\\
Florida Atlantic University, Boca Raton, FL 33431, USA}
\thanks[txw]{E-mail: {\tt ywang@math.fau.edu}}

\maketitle

\begin{abstract}
The fundamental property of strongly monotone systems, and strongly
cooperative systems in particular, is the limit set dichotomy due to Hirsch:
if $x(0) < y(0)$, then either $\omega(x)<\omega (y)$, or $\omega (x)=\omega
(y)$ and both sets consist of equilibria. 
We provide here a counterexample showing that this property need not hold for
(non-strongly) cooperative systems.

\end{abstract}


\end{frontmatter}

\section{Introduction}\label{s-intro}

The field of cooperative, and more generally monotone systems, provides one of
the most fruitful areas of theory as well as practical applications
--particularly in biology-- of dynamical systems.  For an excellent
introduction, see the textbook by Smith~\cite{smith} and the recent
exposition~\cite{Hirsch-Smith}.
One of its central tools is a classical theorem of Hirsch
(\cite{Hirsch1,Hirsch2}), the ``limit set dichotomy'' for strongly monotone
(in particular, strongly cooperative) systems, see Theorem 1.16
in~\cite{Hirsch-Smith}.
The limit set dichotomy states that if $x(0)<y(0)$, then either $\omega
(x)<\omega (y)$, or 
$\omega (x)=\omega (y)$ and both sets consist of equilibria.  

According to the recent survey~\cite{Hirsch-Smith}, the problem of deciding if
there are any cooperative systems for which the dichotomy fails is still open.
In~\cite{Hirsch-Smith}, example 1.24, one finds a system which is
monotone but not strongly monotone, for which the dichotomy fails.
The order in this example is the ``ice cream cone'' order, and the authors
explicitly state that it is unknown whether a polyhedral cone example
exists.  
A cooperative system is one defined by a set of ordinary differential equations
$\dot x=f(x)$, where $f = (f_1,\ldots ,f_n)'$, with the property that
$\frac{\partial f_i}{\partial x_j}(x)\geq 0$ for all $i\not= j$ and all $x$.
Cooperative systems are monotone with respect to a polyhedral cone, namely the
main orthant in $\R^n$.  Thus, a counterexample using cooperative systems
provides an answer to this open question.  We provide such a counterexample
here. 

To be precise, we construct here two differentiable functions
\[
f,g: \R \rightarrow  \R
\]
such that $f(0)=g(0)=0$, $xf(x)<0$ and $yg(y)<0$ for all $x,y\not= 0$, and
consider essentially the following system:
\beqn
\dot x&=&f(x)\\
\dot y&=&g(y)\\
\dot z&=&x+y \,.
\eeqn
This system is cooperative.
Note that solutions of the $x$ and $y$ equations converge to zero as
$t\rightarrow \infty $. 
Moreover, for this system,  there exist $\bar x, \bar y$ 
such that the following property holds: 
\begin{quotation}
There exists some $\delta>0$ such that for any solution
$X$ with initial condition $(x(0),y(0),z(0))^T$ such that
$\abs{x(0)-\bar x}<\delta$,  $\abs{y(0)-\bar y}<\delta $, and
$\abs{z(0)} < 1$,
the omega-limit set $\omega (X)$ 
is compact and it contains the set
\[
 \left\{(0,0,\zeta) \;\bigg|\; 
 z(0)-\frac{1}{2}\leq \zeta\leq z(0)+\frac{1}{2}
 \right\}\,.
 \]
\end{quotation}
The limit set dichotomy states would imply that, for any initial conditions
\[
(x(0),y(0),z(0))^T \quad\mbox{and}\quad
(\hat x(0),\hat y(0),\hat z(0))^T
\]
for which $x(0)\leq \hat x(0)$,
$y(0)\leq \hat y(0)$, and
$z(0)\leq \hat z(0)$,
and with at least one of the inequalities being strict, 
the corresponding solutions $X, \hat X$ have the property that either
$\omega (X)=\omega (\hat X)$ or
$\omega (X)<\omega (\hat X)$.
This last property implies in particular that 
$\omega (X)$ and $\omega (\hat X)$ are disjoint.
Now, for our example, clearly 
$\omega (X)\not= \omega (\hat X)$ as long as $z(0)\not= \hat z(0)$
(since the $z$-components of solutions are translates by $z(0)$
of the solutions with $z(0)=0$ and the $\omega$-limit sets are compact),
but the omega-sets intersect as long as the $x$ and $y$ initial conditions are
less than $\delta $, $\hat z(0)\leq z(0)$, and 
\[
z(0)-\frac{1}{2}
<
\hat z(0) + \frac{1}{2}\,,\quad \abs{z(0)}< 1, \ \abs{\hat z(0)}< 1,
\]
i.e.
\[
0\leq z(0) - \hat z(0) \leq  1\,, \quad \abs{z(0)}< 1,\ \abs{\hat z(0)}< 1.
\]
Thus we contradict the limit set dichotomy.
  
In fact, in the above
discussion, $\bar x$ and $\bar y$ can be chosen arbitrarily small, that
is, for any $r_0>0$, $\bar x$ and $\bar y$ can be chosen so that
$\abs{\bar x} < r_0, \ \abs{\bar y} < r_0$.  Also note that by properly
choosing a function $\sigma(z)$ and modifying  the $z$-subsystem to $\dot z =
x+y-\sigma(z)$, all trajectories of the system can be made to have compact
closures and the above discussions will still remain valid.  For details, see
Proposition \ref{p-main}.

\section{The Example}

To present an example as discussed in Section \ref{s-intro}, we first
consider the following result.

\bl{l-1} For any $\delta>0$, there exist two $C^1$ functions $p,
q: [-1, \infty)\rightarrow [0, \infty)$ 
such that the following holds:
\begin{enumerate}
\item both $p$ and $q$ are strictly decreasing functions;
\item $0< p(0) < \delta$, $0< q(0) < \delta$;
\item $p(t)\rightarrow 0$ and $q(t)\rightarrow 0$;
\item for any $a, b\in (-1, 1)$, the function $H_{a, b}$ defined by 
\[
H_{a, b}(T) = \int_0^T(p(t+a) - q(t+b))\,dt
\]
is bounded; and
\item
for any $a, b\in (-1, 1)$,
 \[
 \mylimsup{T\rightarrow\infty} H_{a, b}(T) >  0,\qquad
 \myliminf{T\rightarrow\infty} H_{a, b}(T) <  0;
\]
and
\[
\mylimsup{T\rightarrow\infty} H_{a, b}(T) -
\myliminf{T\rightarrow\infty} H_{a, b}(T)\ge 1.
\]
\end{enumerate}
\els

We will prove the lemma in Section \ref{s-pf} by showing that the two
functions can be chosen as 
\beq
p(t) &=&\frac{1}{\sqrt{t+c_0}}\label{e-p}\\
q(t) &=&\frac{1}{\sqrt{t+c_0}} +
\frac{\sin[(t+c_0)^{1/4}]}{(t+c_0)^{3/4}},
\label{e-qq}
\eeq
where $c_0$ can be any number in $[\alpha, \infty)$ for some $\alpha$ to
be chosen.  Observe that with proper choices of $c_0$, one can have 
$p(0) = 1/\sqrt{c_0}$, $q(0) = [1/\sqrt{c_0}] + 1/c_0^{3/4}$.
As a consequence, property (ii) in the lemma can be fulfilled with
large enough values of $c_0$.

Below we show that, for some $C^1$ maps $f$ and $g$, 
$p(t+a)$ is a solution of $\dot x = f(x)$; $q(t+a)$ and 
$-q(t+a)$ are solutions of $\dot y = g(y)$ for any $\abs{a}<1$.

It is readily seen that, for $\abs{a} <1$, 
$p(t+a)$ is the solution of the initial
value problem
\[
\dot x = -\frac{x^3}{2}, \quad x(0) = \frac{1}{\sqrt{c_0+a}}.
\]
To get a $C^1$ map $g(\cdot)$ with the desired properties, 
let $\psi: (-1, \infty)\rightarrow (-\infty, 0)$ be defined by 
$\psi(t) = q'(t)$, and
$\vf: (0, \rho)\rightarrow (-1, \infty)$ be
defined by $\vf(r) = q^{-1}(r)$, where $\rho=q(-1)\ge q(t)$ for all
$t > -1$.
Let $g(r) = \psi\circ \vf(r)$ for $r\not= 0$, and $g(0) =0$. 

\bl{l-go}
The function $g: [0, \rho)\rightarrow (-\infty, 0)$ is of $C^1$.
\els

We will prove Lemma \ref{l-go} in Section \ref{s-plgo}.

Extend $g$ from $[0, \rho)$ to $[0, \infty)$ as a $C^1$ function, and
 then extend $g$ to $\R$ by letting $g(-r) = -g(r)$ for $r < 0$.
 Still denote the newly extended function by $g$.  Then $g$ is a $C^1$
 function.  Let $f(x) = -x^3/2$.

\bl{l-2}
For any $a, b\in(-1, 1)$,  the function $(x_a(t), y_b(t))^T$
defined by
\[
x_a(t) = p(t+a)), \quad y_b(t) = -q(t+b) 
\]
is the solution to the initial value problem
\be{e-sysxy}
\begin{array}{lllll}
\dot x &=& f(x), \quad x_a(0) &=& \frac{1}{\sqrt{c_0+a}},\\
\dot y &=& g(y), \quad y_b(0) &=& -q(b).
\end{array}
\ee
\els
The proof of Lemma \ref{l-2} will be given in Section \ref{s-pl2}.

To get a system as discussed in Section \ref{s-intro}, we would like to
cascade the system \rref{e-sysxy} with the one-dimensional system 
$\dot z = x+y$.  To obtain a system for which all trajectories are
bounded, we choose a $C^1$ function $\sigma:\R\rightarrow\R$ with the
property such that 
\begin{itemize}
\item[$\bullet$] $\sigma(r) = 0$ for all $\abs{r}\le 1+M$, where $M =
  \sup\{\abs{H_{a, b}(t)}: \; \abs{a}\le 1, \abs{b}\le 1, t\ge 0\}$; 
\item[$\bullet$] $r\sigma(r) > 0$ for all $\abs{r} > M+1$; and
\item[$\bullet$] $\sigma$ is proper.
\end{itemize}

Consider the system
\be{e-sysxyz}
\begin{array}{lll}
\dot x &=& f(x)\\
\dot y &=& g(y)\\
\dot z &=& x+y - \sigma(z).
\end{array}
\ee

It is clear that the system \rref{e-sysxyz} is cooperative.

\bl{l-z} Consider system \rref{e-sysxyz}:
\begin{enumerate}
\item the $(x, y)$-subsystem is globally asymptotically stable; and
\item every trajectory of the system \rref{e-sysxyz} has a compact closure.
\end{enumerate}
\els

Below we present our final result. We use $X(t) = (x(t), y(t), z(t))^T$
to denote a solution of the system \rref{e-sysxyz}.

\bp{p-main}
Consider the cooperative system \rref{e-sysxyz}.  For any $\delta>0$,
there exist two trajectories $X_1(t)$ and $X_2(t)$ with $X_1(0) < X_2(0)$ and
$\abs{X_1(0)}< \delta, \abs{X_2(0)} < \delta$ such that
$\omega(X_1)\not=\omega(X_2)$ and $\omega(X_1) \le \omega(X_2)$ fails.
\eps

\br{r-1} In fact, we have obtained a system for which the statement of
Proposition \ref{p-main} can be made generic in the following sense.
For any given $\delta_0>0$, there exists $\abs{X_0} < \delta_0$ and some
$\delta_1 > 0$ such that for any pair of trajectories $X_1, X_2$ of the
system \rref{e-sysxyz} satisfying $\abs{X_1(0) - X_0} < \delta_1$,
$\abs{X_2(0) - X_0} < 
\delta_1$, and $z_1(0) \not = z_2(0)$, $\abs{z_1(0)} < 1, \ \abs{z_2(0)}
< 1$, it holds that
$\omega(X_1)\not=\omega(X_2)$ and $\omega(X_1) \le \omega(X_2)$ fails.
\er

{\it Proof of Proposition \ref{p-main}.}\
Assume that $\delta>0$ is given.  Choose $c_0$ as in \rref{e-p}--\rref{e-qq}
large enough so that $\frac{1}{\sqrt{c_0-1}} <\delta$ and $q(-1) <
\delta$.  According to Lemma \ref{l-2}, 
for any trajectory of the system
\[
\dot x = f(x), \quad \dot y =  g(y)
\]
with 
\be{e-xoyo}
\frac{1}{\sqrt{c_0+1}}< x(0) < \frac{1}{\sqrt{c_0-1}},\qquad
-q(-1) < y(0) < -q(1),
\ee
one has $x(t) = p(t+a)$, $y(t) = -q(t+b)$ for some $a, b\in(-1, 1)$, and
hence, 
\be{e-hpq}
\begin{array}{l}
{\displaystyle \mylimsup{t\rightarrow\infty}\int_0^t (x(s) + y(s))\,ds - 
\myliminf{t\rightarrow\infty}\int_0^t (x(s) + y(s))\,ds}\\
=
\mylimsup{t\rightarrow\infty} H_{a, b}(t) -
\myliminf{t\rightarrow\infty} H_{a, b}(t) \ge 1.
\end{array}
\ee
Let $x_0 = 1/\sqrt{c_0},\  y_0 = q(0)$, and 
\[
\delta_1 = \min\left\{\frac{1}{\sqrt{c_0-1}} - \frac{1}{\sqrt{c_0}},\;
\frac{1}{\sqrt{c_0}} - \frac{1}{\sqrt{c_0+1}},\;
q(-1) - q(0), \; q(0) - q(1)\right\}.
\]
Then \rref{e-hpq} holds for any trajectory $(x(t), y(t))^T$ of
\rref{e-xoyo} with $\abs{x(0)-x_0}<\delta_1, \ \abs{y(0) - y_0} <
\delta_1$.
Take any $\abs{z_0} < 1$, 
\[
\abs{z_0 + \int_0^t (x(s) + y(s))\, ds}
= \abs{z_0 + H_{a, b}(t)}\ < M +1,
\]
and therefore, $z_0 + \int_0^t (x(s) + y(s))\, ds$ is the solution of
the $z$-subsystem of \rref{e-sysxyz} with the initial value $z(0) = z_0$
(note that $\sigma(s) = 0$ when $\abs{s}\le M+1$).  It then follows from
statement (v) of Lemma \ref{l-1} that for any $x(0), y(0)$
satisfying \rref{e-xoyo}, it holds that
\[
\mylimsup{t\rightarrow\infty}\,(z(t) - z(0)) > 0, \quad
\myliminf{t\rightarrow\infty}\,(z(t) - z(0)) < 0,
\]
and
\[
\left[\mylimsup{t\rightarrow\infty}\, (z(t) - z(0))\right]
-
\left[\myliminf{t\rightarrow\infty}\, (z(t) - z(0))\right] > 1.
\]
Observe that for any trajectory $X(t)$ of the system \rref{e-sysxyz},
\[
\omega(X) \supseteq \left\{(0, 0, \alpha):\; 
\myliminf{t\rightarrow\infty}\,z(t) < \alpha
<\mylimsup{t\rightarrow\infty}\,z(t)\right\}.
\]
For any two solutions
$X_1(t) :=(x(t), y(t), \hat z(t))$ and $X_2(t) : =(x(t), y(t), \hat
z(t))$, where $x(0)$ and $y(0)$ satisfy \rref{e-xoyo}, and $\hat
z(0)\not= z(0)$, $\abs{z(0)} < 1, \ \abs{\hat z(0)} <1$, 
the sets of $\omega$-limit
  points of $z(t)$ and $\hat z(t)$ are different
  (since $z(t) - \hat z(t) \equiv z(0) - \hat z(0)$).  
Moreover, with  $\hat z(0) < z(0)$,
\beqn
\mylimsup{t\rightarrow\infty}\, \hat z(t) -
\myliminf{t\rightarrow\infty} z(t)& = &
\left[\mylimsup{t\rightarrow\infty}\, (\hat z(t) - \hat z(0))\right]
-
\left[\myliminf{t\rightarrow\infty}\, (z(t) - z(0))\right]\\
&-& (z(0) - \hat z(0))
\ge 1 - (z(0) - \hat z(0)) > 0
\eeqn
if $0< z(0) - \hat z(0) < 1$.  Hence, $\omega(\hat X) \le \omega( X)$
fails. 

\section{Proofs of the Lemmas}
In this section, we provide proofs of the results.

\subsection{Proof of \protect Lemma \ref{l-1}}\label{s-pf}

First we let
\be{e-qo}
p_0(t) = \frac{1}{\sqrt{t}},\quad
q_0(t) = \frac{1}{\sqrt{t}} + \frac{1}{t^{3/4}}\sin t^{1/4}\quad t\ge 1.
\ee
To show that $q_0$ is decreasing, we consider $q_0'(t)$:
\beq
q_0'(t) &=& -\frac{1}{2t^{3/2}} - \frac{3}{4t^{7/4}}\sin t^{1/4} + 
\frac{1}{t^{3/4}}\cdot\frac{1}{4 t^{3/4}}\cos t^{1/4}\nonumber\\[3mm]
& =&
-\frac{1}{2t^{3/2}} - \frac{3}{4t^{7/4}}\sin t^{1/4} + 
\frac{1}{4t^{3/2}}\cos t^{1/4}\nonumber\\[3mm]
&\le&
-\frac{1}{4t^{3/2}} + \frac{3}{4t^{7/4}}\label{e-g}.
\eeq
So, $q_0'(t) \le 0$ when $\dis{\frac{1}{t^{3/2}} - \frac{3}{t^{7/4}}\ge
  0}$.  This is the same as $\dis{t^{3/2}\le t^{7/4}/3}$, or $t^{1/4} \ge
3$, that is, $t\ge 81$.

Let $p(t) = p_0(t+c_0)$, $q(t) = q_0(t+c_0)$, where 
$c_0\ge 81$ will be chosen later on. Now both $p$ and $q$ are
differentiable on $[0, \infty)$ and monotonically decrease to $0$.
For $a\in (-1, 1)$ and $b\in (-1, 1)$,  let 
$H_{a, b}(T)$ be as defined as in Lemma \ref{l-1}. Then
\beq
& &H_{a, b}(T)
= \int_0^T (p(t+a) -q(t+b)\,dt\nonumber\\
&=&\int_0^T \frac{1}{\sqrt{t+c_0+a}} - \frac{1}{\sqrt{t+c_0+b}}
-\frac{1}{(t+c_0+b)^{3/4}}\sin(t+c_0+b)^{1/4}\,dt\nonumber\\[3mm]
&=&
\int_0^T\frac{\sqrt{t+c_0+b}-
  \sqrt{t+c_0+a}}{\sqrt{t+c_0+a}\,\sqrt{t+c_0+b}} 
-\frac{1}{(t+c_0+b)^{3/4}}\sin(t+c_0+b)^{1/4}\,dt.\label{e-Hab}
\eeq
Now, for the first term in \rref{e-Hab}, we have
\beqn
& &\left|\frac{\sqrt{t+c_0+b}-
  \sqrt{t+c_0+a}}{\sqrt{t+c_0+a}\sqrt{t+c_0+b}} \right|\\
&=& \frac{b-a}{{\sqrt{t+c_0+a}\sqrt{t+c_0+b}}(\sqrt{t+c_0+b}+
  \sqrt{t+c_0+a})}\\
&\le&\frac{|b-a|}{2(t+c_0-1)^{3/2}} \qquad\quad\forall \abs{a}<1,
\abs{b}<1, 
\eeqn
and hence, the integral 
\[
\int_0^\infty\frac{\sqrt{t+c_0+b}-
  \sqrt{t+c_0+a}}{\sqrt{t+c_0+a}\sqrt{t+c_0+b}}\,dt
\]
is convergent, and for $\abs{a} < 1, \abs{b}<1$,
\be{e-H1}
\abs{\int_0^T \frac{\sqrt{t+c_0+b}-
  \sqrt{t+c_0+a}}{\sqrt{t+c_0+a}\sqrt{t+c_0+b}}\,dt}\le
\frac{2}{\sqrt{c_0}}\quad\forall\,T > 0.
\ee
For the second term in \rref{e-Hab}, using
$u=(t+c_0+b)^{1/4}$, one has
\beqn
& &\int_0^T \frac{1}{(t+c_0+b)^{3/4}}\sin(t+c_0+b)^{1/4}\,dt
= \int_{(c_0+b)^{1/4}}^{(T+c_0+b)^{1/4}}\sin u\,du\\
& & = -\cos(T+c_0+b)^{1/4} +
\cos (c_0+b)^{1/4}.
\eeqn
Combining this with \rref{e-H1}, one sees that, for any $\abs{a}<1,
\abs{b}<1$, $H_{a, b}(t)$ is bounded on $[0, \infty)$.

Let $c_0\ge 82$ be of the form $(2k\pi+\pi/2)^4$.  For
$\abs{b}< 1 $,
\[
\frac{d}{db}(c_0+b)^{1/4} = \frac{1}{4(c_0+b)^{3/4}}
\le\frac{1}{4}\le \frac{\pi}{6}
\]
Consequently,
\[
c_0^{1/4} - \frac{\pi}{6}
\le (c_0+b)^{1/4}\le c_0^{1/4} + \frac{\pi}{6}\qquad\forall\,\abs{b} < 1.
\]
This implies that
\[
-\frac{1}{2}\le
\cos(c_0+b)^{1/4}\le\frac{1}{2}\qquad\forall\,\abs{b} < 1.
\]
Thus,
\beqn
& &\mylimsup{T\rightarrow\infty}
\int_0^T \frac{1}{(t+c_0+b)^{3/4}}\sin(t+c_0+b)^{1/4}\,dt =
1 - \cos(c_0+b)^{1/4}\ge \frac{1}{2}\\
& &
\myliminf{T\rightarrow\infty}
\int_0^T \frac{1}{(t+c_0+b)^{3/4}}\sin(t+c_0+b)^{1/4}\,dt 
= -1 - \cos(c_0+b)^{1/4}\le - \frac{1}{2}.
\eeqn

Finally, we let $c_0 = (2k\pi + \pi/2)^4$ with $k$ large enough so that
$c_0 \ge 82$ (and consequently $\frac{2}{\sqrt{c_0}} < 1/4$).
This way, we get for all $\abs{a}<1, \abs{b}< 1$,
\beqn
\mylimsup{T\rightarrow\infty} H_{a, b}(T)&=& \lim_{T\rightarrow\infty}
\int_0^T\left(\frac{1}{\sqrt{t+c_0+a}} -
  \frac{1}{\sqrt{t+c_0+b}}\right)\,dt \\
&-&\myliminf{T\rightarrow\infty}
\int_0^T\frac{1}{(t+c_0+b)^{3/4}}\sin(t+c_0+b)^{1/4}\,dt\\
&\ge& l_{a, b} + \frac{1}{2}, 
\eeqn
where $\dis{l_{a, b} = \int_0^\infty\left(\frac{1}{\sqrt{t+c_0+a}}
    -\frac{1}{\sqrt{t+c_0+b}}\right)\, dt}$; and
\beqn
\myliminf{T\rightarrow\infty}H_{a, b}(T)&= & 
\lim_{T\rightarrow\infty}
\int_0^T\left(\frac{1}{\sqrt{t+c_0+a}} -
\frac{1}{\sqrt{t+c_0+b}}\right)dt \\ 
& - & \mylimsup{T\rightarrow}
 \int_0^T\frac{1}{(t+c_0+b)^{3/4}}\sin(t+c_0+b)^{1/4}\,dt\\
&\le& l_{a, b} - \frac{1}{2}, 
\eeqn
which implies that
\[
\mylimsup{t\rightarrow\infty}H_{a, b}(t) -
\myliminf{t\rightarrow\infty}H_{a, b}(t) \ge 1.
\]
Since $c_0$ was chosen so that $\abs{l_{a, b} } < 1/4$, one has
\[
\myliminf{T\rightarrow\infty}H_{a, b}(T) < -1/4,\quad
\myliminf{T\rightarrow\infty}H_{a, b}(T) > 1/4.
\] 

\subsection{Proof of \protect{Lemma \ref{l-go}}}\label{s-plgo}

First of all, it can be calculated (see also
\rref{e-g}) that,  for $t$ large enough
\be{e-psi}
\abs{\psi(t)} =\abs{\frac{d}{dt}q_0(t+c_0)} \le \frac{M}{t^{3/2}}
\ee
for some $M\ge 0$, where $q_0$ is defined as in \rref{e-qo}.
Let $\vf: (0, \rho)\rightarrow (-1, \infty)$ be
defined by $\vf(r) = q^{-1}(r)$, where $\rho=q(-1)\ge q(t)$ for all
$t > -1$.  Note that
\[
\lim_{r\rightarrow 0}\vf(r) = \infty.
\]
Let $g(r) = \psi\circ \vf(r)$ for $r\not= 0$, and $g(0) =0$.  Then $g$
is continuous on $[0, \rho)$, and of $C^1$ on $(0, \rho)$. Observe that
$g(r) < 0$ for all $r\in (0, \rho)$. Below we
show that $g$ is differentiable at $0$. 

\noindent{\it Fact 1.} There exist some $\delta_0>0$ and some $L_0>0$ 
such that 
\be{e-vf1}
\vf(r) \ge \frac{L_0}{r^2}\qquad\forall\,r\in (0, \delta_0).
\ee

To prove Fact 1, write $t = q^{-1}(r)$. Then $r = q(t)$, and
\be{e-r}
r \ge \frac{1}{\sqrt{t+c_0}}  - \frac{1}{(t+c_0)^{3/4}}\ge
\frac{L_0}{\sqrt{t}} \qquad \forall\, t\ge T_0
\ee
for some $L_0>0$ and some $T_0\ge 0$.  Since $q(t)$ decreases to $0$ as
$t\rightarrow\infty$, it follows that for some $\delta_0>0$, it holds
that $t\ge \frac{L_0}{r^2}$ for all $r\in (0, \delta_0)$, this is,
\rref{e-vf1} holds.

\noindent{\it Fact 2.} $g'(0) = 0$.

Fact 2 follows from Fact 1 combined with \rref{e-psi}:
\[
\abs{\psi(\vf(r))}\le \frac{M}{[\vf(r)]^{3/2}}\le
M\left[\frac{L_0}{r^2}\right]^{-3/2} \le \tilde Mr^3
\]
for all $r>0$ in a neighborhood of $0$, where $\tilde M>0$ is some
constant. This shows that $g$ is differentiable at $0$, and $g'(0) = 0$. 

\noindent{\it Fact 3.} $q'(r)$ is continuous at $r=0$.

To prove this fact, we first get an estimate on $\psi'(t)$ for $t$ large
enough:
\beqn
\psi'(t-c_0) = q_0''(t) &=& \frac{3}{4t^{5/2}} + \frac{21}{16 t^{11/4}}\sin
  t^{1/4} \\
&-& \frac{3}{16t^{10/4}}\cos t^{1/4} - \frac{3}{8t^{5/2}}\cos
    t^{1/4}  - \frac{1}{16t^{9/4}}\sin t^{1/4}.
\eeqn
Hence, for some $T_1>0$ and some $L_1\ge 0$,
\[
\abs{\frac{d}{dt}\psi(t)}\le \frac{L_1}{t^{9/4}}\qquad \forall\,t\ge T_1.
\]
This implies that, for some $\delta_1>0$,
\be{e-q}
\abs{\psi'(\vf(r))}\le \frac{L_1}{[\vf(r)]^{9/4}}\qquad\forall\,r\in (0,
  \delta_1). 
\ee
We also need the following estimate on $\vf'(r)$ near $0$:
\beqn
\abs{\frac{d}{dt}q(t-c_0)} &=& \abs{\frac{d}{dt}q_0(t)}\ge 
\frac{1}{2t^{3/2}}
-\frac{3}{4t^{7/4}} - \frac{1}{4t^{3/2}}\\
&=& \frac{1}{4t^{3/2}} - \frac{3}{4t^{7/4}}\qquad\qquad\forall\, t\ge c_0. 
\eeqn
It then can be seen that for some $T_2>0$ and some $L_2 >0$, one has
\[
\abs{q'(t)}\ge \frac{L_2}{t^{3/2}}\qquad\forall\, t\ge T_2,
\]
which implies that for some $\delta_2 >0$,
\[
\abs{q'(\vf(r))}\ge \frac{L_2}{[\vf(r)]^{3/2}}\qquad\forall\,r\in(0,
\delta_2). 
\]
Finally,
\beqn
\abs{g'(r)} &=&\abs{\psi'(\vf(r))\vf'(r)} = 
\abs{\psi'(\vf(r))\frac{1}{q'(\vf(r))}}\\[2mm]
&\le&
\frac{L_1}{[\vf(r)]^{9/4}}\cdot \frac{1}{\frac{L_2}{[\vf(r)]^{3/2}}}
=
\frac{L_1}{L_2[\vf(r)]^{3/4}}\rightarrow 0\quad {\rm as}\
r\rightarrow 0.
\eeqn
Hence, $g'(r)$ is continuous at $r=0$.  With this we conclude that $g$
is of $C^1$ on $[0, \rho)$.~\qed

\subsection{Proof of \protect{Lemma \ref{l-2}}}\label{s-pl2}

The statement about $x_a(t)$ is certainly clear.  

To treat the part about $y_b(t)$, first 
observe that $q'(t)$ can be written as $q'(\vf(q(t)))$.  Also note that
for any $\abs{b}<1$, $0\le q(t+b)\le q(-1)$ for all $t\ge 0$, that
is, $q(t+b)\in (0, \rho)$ for all $t\ge 0$. Hence, we have
\[
\frac{d}{dt}q(t) = g(q(t)),
\]
that is, $q(t)$ is a solution of the differential equation $\dot y= q(y)$.
Let $\tilde q(t) = -q(t)$. Then
\[
\frac{d \tilde q(t)}{dt}  = -\frac{dq}{dt} = -g(q(t)) = -g(-\tilde q(t))
= g(\tilde q(t)).
\]
This shows that $-q(t)$ is also a solution of the equation
$\dot y= g(y)$.~\qed

\subsection{Proof of \protect{Lemma \ref{l-z}}}

Since for $xf(x) <0$ and $yg(y)< 0$ for all $x\not=0, y\not= 0$, both the
$x$- and the $y$-subsystems are globally asymptotically stable.

To complete the proof, it is enough to note that every trajectory of the
system 
\[
\dot z = -\sigma(z) + h(t)
\]
is bounded for any choice of bounded function $h(t)$.~\qed

\end{document}